\documentclass[12pt,a4paper]{article} 

\usepackage{amsmath, amssymb, verbatim} 
\usepackage{theorem} 
 
\newtheorem{theorem}{Theorem}[section]

\newtheorem{proposition}[theorem]{Proposition}
\theorembodyfont{\normalfont} 
\newtheorem{definition}[theorem]{Definition}
\newtheorem{remark}[theorem]{Remark}
\newtheorem{example}[theorem]{Example}
\newtheorem{conjecture}[theorem]{Conjecture}

\newcommand{\proj}{\mathbb{P}} 
\newcommand{\cpx}{\mathbb{C}} 
\newcommand{\om}{\omega} 
\newcommand{\eps}{\varepsilon}

\newcommand{\ra}{\rightarrow} 
\newcommand{\olo}{\mathcal{O}}

\newcommand{\Bl}{\textrm{Bl}}
\newcommand{\Sl}{\textrm{SL}}
\numberwithin{equation}{section}   
 
\newcommand{\Aut}{\operatorname{Aut}} 

\newcommand{\dimo}[1][]         {\noindent\textbf{Proof#1}. } 

\newcommand{\fine}    {\begin{flushright} 
             \textsc{Q.E.D.} 
             \end{flushright}}

\begin{document} 
\title{K-stability of constant scalar curvature K\"ahler manifolds}

\author{Jacopo Stoppa} 
\date{} 
 
\maketitle 
 
\begin{abstract} 
  \noindent 
We show that a polarised manifold with a constant scalar curvature K\"ahler metric and discrete automorphisms is K-stable. This refines the K-semistability proved by S. K. Donaldson.
\end{abstract} 
\section{Introduction}
Let $(X, L)$ be a polarised manifold. One of the more striking realisations in K\"ahler geometry over the past few years is that if one can find a constant scalar curvature K\"ahler (cscK) metric $g$ on $X$ whose $(1,1)$-form $\om_g$ belongs to the cohomology class $c_1(L)$ then $(X, L)$ is \textit{semistable}, in a number of senses. The seminal references are Yau \cite{yau}, Tian \cite{tian}, Donaldson \cite{don_fields}, \cite{don_toric}. 

In this note we are concerned with Donaldson's \textit{algebraic K-stability} \cite{don_toric}, see also Definition \ref{semi} below.  This notion generalises Tian's K-stability for Fano manifolds \cite{tian}. It should play a role similar to Mumford-Takemoto slope stability for bundles. The necessary general theory is recalled in Section \ref{gen_theory}.

Asymptotic Chow stability (which implies K-semistability, see e.g. \cite{ross_thomas} Theorem 3.9) for a cscK polarised manifold was first proved by Donaldson \cite{don_scalar} in the absence of continous automorphisms. Important work in this connection was also done by Mabuchi, see e.g. \cite{mabuchi}. From the analytic point of view the fundamental result is the lower bound on the K-energy proved by Chen-Tian \cite{chen_tian}.

The neatest result in the algebraic context seems to be Donaldson's \emph{lower bound on the Calabi functional}, which we now recall.

For a K\"ahler form $\om$ let $S(\om)$ denote the scalar curvature, $\widehat{S}$ its average (a topological quantity). Denote by $F$ the Donaldson-Futaki invariant of a test configuration (Definitions \ref{test_config}, \ref{futaki}). The precise definition of the norm $\|\mathcal{X}\|$ appearing below will not be important for us. 
\begin{theorem}[Donaldson \cite{don_calabi}]\label{calabi} For a polarised manifold $(X, L)$
\begin{equation}\label{calabi_ineq} 
\inf_{\om \in c_1(L)}\int_X (S(\om)-\widehat{S})^2\om^n \geq -\frac{\sup_{\mathcal{X}} F(\mathcal{X})}{\|\mathcal{X}\|}.
\end{equation}
where the supremum is taken with respect to all test configurations $(\mathcal{X}, \mathcal{L})$ for $(X, L)$. 

Thus if $c_1(L)$ admits a cscK representative $(X, L)$ is K-semistable.
\end{theorem}
There is a strong analogy here with Hermitian Yang-Mills metrics on vector bundles. By the celebrated results of Donaldson and Uhlenbeck-Yau these are known to exist if and only if the bundle is slope polystable, namely a semistable direct sum of slope stable vector bundles. 

In particular a simple vector bundle endowed with a HYM metric is slope stable. In this note we will prove the corresponding result for polarised manifolds. 
\begin{theorem}\label{main_teo} If $c_1(L)$ contains a cscK metric and $\Aut(X, L)$ is discrete then $(X, L)$ is K-stable. 
\end{theorem}
Theorem \ref{main_teo} fits in a more general well known conjecture.
\begin{conjecture}[Donaldson \cite{don_toric}]\label{don_conj} If $c_1(L)$ contains a cscK metric then $(X, L)$ is K-polystable (Definition \ref{poly}).  
\end{conjecture}
Thus our result confirms this expectation when the group $\Aut(X,L)$ is discrete. From a differential-geometric point of view this means that $X$ has no nontrivial Hamiltonian holomorphic vector fields - holomorphic fields that vanish somewhere.
\begin{remark} Conjecture \ref{don_conj} and its converse are known as Yau - Tian - Donaldson Conjecture, and sometimes called the Hitchin-Kobayashi correspondence for manifolds.
\end{remark} 
For the rest of the note we will assume $\dim(X) > 1$ in all our statements.

K-stability for Riemann surfaces is completely understood thanks to the work of Ross and Thomas \cite{ross_thomas} Section 6. In particular Conjecture \ref{don_conj} is known to hold for Riemann surfaces.\\ 

Our proof of Theorem \ref{main_teo} rests on the general principle that one should be able to \emph{perturb} a semistable object (in the sense of geometric invariant theory) to make it unstable - altough this necessarily involves perturbing the GIT problem too, since the locus of semistable points for an action on a fixed variety is open. Conversely in the absence of continous automorphisms, the cscK property is open - at least in the sense of small deformations - so cscK should imply stability. Of course we need to make this rigorous; in particular testing small deformations is not enough to prove K-stability.\\

Thus suppose that $(X, L)$ is properly K-semistable (Definition \ref{prop}). We will find a natural way to construct from this a family of K-unstable small perturbations $(X_{\eps}, L_{\eps})$ for small $\eps > 0$. Our choice for $X_{\eps}$ is actually constant, the blowup $\widehat{X} = \Bl_q X$ at a very special point $q$ with exceptional divisor $E$. Only the polarisation changes, and quite naturally $L_{\eps} = \pi^*L - \eps\olo(E)$. This would involve taking $\eps \in \mathbb{Q}^{+}$ and working with $\mathbb{Q}$-divisors, but in fact we rather take tensor powers and work with $\widehat{X}$ polarised by $L_{\gamma} = \pi^* L^{\gamma} - \olo(E)$ for integer $\gamma \gg 0$. K-(semi, poly, in)stability is unaffected by Definition \ref{futaki}. 
\begin{proposition}\label{defo} Let $(X, L)$ be a properly K-semistable polarised manifold. Then there exists a point $q \in X$ such that the polarised blowup $(\emph{Bl}_q X, \pi^*L^{\gamma}\otimes\olo(-E))$ is K-unstable for $\gamma \gg 0$.
\end{proposition}
\begin{remark} It is interesting to note that the corresponding result for vector bundles follows from Buchdahl \cite{buch}. Let $(X, L)$ be a polarised manifold and $E \ra X$ a properly slope semistable vector bundle. Then the pullback $\pi^*E$ to the blowup $\Bl_{q_1, ...\,q_m}X$ in a finite number of suitably chosen points is slope unstable with respect to the polarisation $\pi^*L^{\gamma}\otimes\olo_{\Bl_{\{q_i\}}X}(1)$ for $\gamma \gg 0$.
\end{remark}
Assume now that a properly semistable $(X, L)$ also admits a cscK metric $\om \in c_1(L)$. If $\Aut(X,L)$ is discrete the blowup perturbation problem for $\om$ is unobstructed by a theorem of Arezzo and Pacard \cite{arezzoI}, so we would get cscK metrics in $c_1(\pi^*L^{\gamma}\otimes\olo(-E))$ for $\gamma \gg 0$, a contradiction. 
\begin{remark}
This perturbation strategy for proving \ref{don_conj} is very general, and was first pointed out to the author by S. Donaldson and G. Sz\'ekelyhidi. Different choices for $(X_{\eps}, L_{\eps})$ lead to different perturbation problems for $\om$, which may settle Conjecture \ref{don_conj} in the presence of continous automorphisms. A possible variant is to perturb the cscK equation with $\eps$ at the same time, but one would then need to develop the relevant K-stability theory for a more general equation. 
\end{remark}
To sum up the main ingredients for our proof (besides Theorem \ref{calabi}) are:
\begin{enumerate}
\item A well known embedding result for test configurations (Proposition \ref{embed}), together with the algebro-geometric estimate Proposition \ref{repulsive};
\item A blowup formula for the Donaldson-Futaki invariant proved by the author \cite{io} Theorem 1.3;
\item A special case of the results of Arezzo and Pacard on blowing up cscK metrics \cite{arezzoI}. 
\end{enumerate}
\textbf{Aknowledgements.} I am grateful to S. K. Donaldson, G. Sz\'ekelyhidi and my advisor R. Thomas for many useful discussions. The reference \cite{buch} was pointed out to me by J. Keller.  
\section{Some general theory}\label{gen_theory}
Let $n$ denote the complex dimension of $X$.
\begin{definition}[Test configuration.]\label{test_config} A \emph{test configuration} for a polarised manifold $(X, L)$ is a polarised flat family $(\mathcal{X}, \mathcal{L})\ra \cpx $ with $(\mathcal{X}_1, \mathcal{L}_1) \cong (X, L)$ and which is  $\cpx^*$-equivariant with respect to the natural action of $\cpx^*$ on $\cpx$.
\end{definition}
Given a test configuration $(\mathcal{X}, \mathcal{L})$ for $(X, L)$ denote by $A_k$ the matrix representation of the induced $\cpx^*$-action on $H^0(\mathcal{X}_0, \mathcal{L}^k_0)$. By (equivariant) Riemann-Roch we can find expansions
\begin{align}
h^0(\mathcal{X}_0, \mathcal{L}^k_0) &= a_0 k^n + a_1 k^{n-1} + O(k^{n-2}),\\
\textrm{tr}(A_k) &= b_0 k^{n+1} + b_1 k^{n} + O(k^{n-1}).
\end{align}
\begin{definition}[Donaldson-Futaki invariant.]\label{futaki} This is the rational number
\begin{equation}
F(\mathcal{X}) = a^{-2}_0(b_0 a_1 - a_0 b_1)
\end{equation}
which is independent of the choice of a lifting of the action to $\mathcal{L}_0$. 

Equivalently $F(\mathcal{X})$ is the coefficient of $k^{-1}$ in the Laurent series expansion of the quotient 
\[\frac{\textrm{tr}(A_k)}{k h^0(\mathcal{X}_0, \mathcal{L}^k_0)}.\]
Note moreover that $F$ is invariant under taking tensor powers, i.e.
\[F(\mathcal{X}, \mathcal{L}) = F(\mathcal{X}, \mathcal{L}^r).\]
Therefore for the rest of this note we will assume without loss of generality that $\mathcal{L}$ is \emph{relatively very ample}.
\end{definition}
\begin{remark}[Coverings]\label{covering} Given a test configuration $(\mathcal{X}, \mathcal{L})$ we can construct a new test configuration for $(X, L)$ by pulling $\mathcal{X}$ and $\mathcal{L}$ back under the $d$-fold ramified covering of $\cpx$ given by $z \mapsto z^d$. This changes $A_k$ to $d\cdot A_k$ and consequently $F$ to $d\cdot F$. 
\end{remark}
\begin{definition} A test configuration $(\mathcal{X}, \mathcal{L})$ is called a \emph{product} if it is $\cpx^*$-equivariantly isomorphic to the product $(X\times\cpx, p^*_X L)$ endowed with the composition of a $\cpx^*$-action on $(X, L)$ with the natural action of $\cpx^*$ on $\cpx$.

A product test configuration is called \emph{trivial} if the associated action on $(X, L)$ is trivial.
\end{definition}
The Donaldson-Futaki invariant $F(\mathcal{X})$ in this case coincides with the classical Futaki invariant for holomorphic vector fields.
\begin{definition}[K-stability]\label{semi} A polarised manifold $(X, L)$ is \emph{K-semistable} if for all test configurations $(\mathcal{X}, \mathcal{L})$ 
\[F(\mathcal{X}) \geq 0.\]
It is \emph{K-stable} if the strict inequality holds for nontrivial test configurations.
\end{definition}
In particular if $(X, L)$ is K-stable $\Aut(X, L)$ must be discrete. The correct notion to take care of continous automorphisms is K-polystability.
\begin{definition}\label{poly}
A polarised manifold $(X, L)$ is \emph{K-polystable} if it is K-semista-\\ble and moreover any test configuration $(\mathcal{X}, \mathcal{L})$ with $F(\mathcal{X}) = 0$ is a product.
\end{definition}
\begin{definition}\label{prop} A polarised manifold $(X,L)$ is \emph{properly K-semistable} if it is K-semistable and it admits a nonproduct test configuration with vanishing Donaldson-Futaki invariant. 
\end{definition}
\begin{remark} The terminology \emph{strictly K-semistable} is also found in the literature with the same meaning.
\end{remark}
Test configurations are well known to be equivalent to 1-parameter flat families induced by projective embeddings.
\begin{proposition}[see e.g. Ross-Thomas \cite{ross_thomas} 3.7]\label{embed} A test configuration for $(X, L)$ is equivalent to a 1-parameter subgroup of $\emph{GL}(H^0(X, L)^*)$.
\end{proposition}
In \cite{io} the author proved a blowup formula for the Donaldson-Futaki invariant. The statement involves some more terminology.
\begin{definition}[Hilbert-Mumford weight.] Let $\alpha$ be a 1-parameter subgroup of $\textrm{SL}(N+1)$, inducing a $\cpx^*$-action on $\proj^N$. Choose projective coordinates $[x_0: ... :x_N]$ such that $\alpha$ is given by $\textrm{Diag}(\lambda^{m_0}, ... \lambda^{m_N})$. The Hilbert-Mumford weight of a closed point $q \in \proj^N$ is defined by
\[\mu(q, \alpha) = -\min\{m_i: q_i \neq 0\}.\] 
Note that this coincides with the weight of the induced action on the fibre of the hyperplane line bundle $\olo(1)$ over the specialisation $\lim_{\lambda \ra 0}\lambda\cdot q$.
\end{definition}
\begin{definition}[Chow weight.]\label{chow} Let $(Y, L)$ be a polarised scheme, $y \in Y$ a closed point, and $\alpha$ a $\cpx^*$-action on $(Y, L)$. Suppose that $L$ is very ample and $\alpha \hookrightarrow \Sl(H^0(Y, L)^*)$. The Chow weight $\mathcal{CH}_{(Y, L)}(q,\alpha)$ is defined to be the Hilbert-Mumford weight of $y \in \proj(H^0(Y, L)^*)$ with respect to the induced action. The definition extends to $0$-dimensional cycles on $Y$, that is effective linear combinations of closed points.
\end{definition}
\begin{theorem}[S. \cite{io} 1.3]\label{blowup_auto} For points $q_i \in X$ and integers $a_i > 0$ let $Z \subset X$ be the 0-dimensional closed subscheme $Z = \cup_i a_i q_i$. Let $\Lambda$ be the 0-cycle on $X$ given by $\sum_i a^{n-1}_i q_i$.

A 1-parameter subgroup $\alpha\hookrightarrow\Aut(X, L)$ induces a test configuration $(\widehat{\mathcal{X}}, \widehat{\mathcal{L}})$ for $(\emph{Bl}_{Z}X, \pi^*L^{\gamma}\otimes\olo_{\emph{Bl}_{Z}X}(1) )$, where $\olo_{\emph{Bl}_{Z}X}(1)$ denotes the exceptional invertible sheaf. More precisely let $O(Z)^-$ be the closure of the orbit of $Z$. Then $\widehat{\mathcal{X}} = \emph{Bl}_{O(Z)^-}\mathcal{X}$ and $\widehat{\mathcal{L}} = \pi^*\mathcal{L}^{\gamma}\otimes\olo_{\widehat{\mathcal{X}}}(1)$.

Suppose that $\alpha$ acts through $\emph{SL}(H^0(X, L)^*)$ with Futaki invariant $F(X)$. Then the following expansion holds as $\gamma \ra \infty$
\[F(\widehat{\mathcal{X}}) = F(X)-\mathcal{CH}_{(X, L)}(\Lambda,\alpha)\frac{\gamma^{1-n}}{2(n-1)!} + O(\gamma^{-n}).\]
\end{theorem}
We will need a slight generalisation of this result, covering blowups of non-product test configurations.
\begin{proposition}\label{blowup} Let $(\mathcal{X}, \mathcal{L})$ be a test configuration for $(X, L)$, $Z = \cup_i a_i q_i$ as above. There is a test configuration $(\widehat{\mathcal{X}}, \widehat{\mathcal{L}})$ for $(\emph{Bl}_{Z}X, \pi^*L^{\gamma}\otimes\olo_{\emph{Bl}_{Z}X}(1))$ with total space $\widehat{\mathcal{X}}$ given by the blowup of $\mathcal{X}$ along $O(Z)^-$. The linearisation is the natural one induced on $\widehat{\mathcal{L}} = \pi^*\mathcal{L}^{\gamma}\otimes\olo_{\widehat{\mathcal{X}}}(1)$. 

Let $q_{i,0} = \lim_{\lambda \ra 0} \lambda\cdot q_i$ be the specialisation, $\Lambda_0$ the $0$-cycle on $\mathcal{X}_0$ given by $\sum_i a^{n-1}_i q_{i,0}$. 

Let $\alpha$ denote the induced action on $(\mathcal{X}_0, \mathcal{L}_0)$ and suppose that $\alpha$ acts through $\emph{SL}(H^0(\mathcal{X}_0, \mathcal{L}_0)^*)$. Then the expansion
\[F(\widehat{\mathcal{X}}) = F(\mathcal{X})-\mathcal{CH}_{(\mathcal{X}_0, \mathcal{L}_0)}(\Lambda_0,\alpha)\frac{\gamma^{1-n}}{2(n-1)!} + O(\gamma^{-n})\]
holds as $\gamma \ra \infty$.  
\end{proposition}
We emphasise that the relevant Chow weight is computed \emph{on the central fibre} $(\mathcal{X}_0, \mathcal{L}_0)$ with its induced $\cpx^*$-action.\\
\dimo The argument of \cite{io} Section 4 goes over verbatim to non-product test configurations, with only two exceptions:
\begin{enumerate}
\item The proof of flatness of the composition $\widehat{\mathcal{X}}\ra\mathcal{X}\ra\cpx$;
\item The identification of the weight $\mathcal{CH}_{(\mathcal{X}_0, \mathcal{L}_0)}(\Lambda_0)$ (with respect to the induced action on $\mathcal{X}_0$) with $\mathcal{CH}_{(X, L)}(\Lambda,\alpha)$.
\end{enumerate}
We do not need the latter identification, and indeed it does not make sense in this case since the general fibre is not preserved by the $\cpx^*$-action.

To prove flatness we use the criterion \cite{hart} III Proposition 9.7. Thus we need to prove that all associated points of $\widehat{\mathcal{X}}$ (i. e. irreducible components and their thickenings) map to the generic point of $\textrm{Spec}(\cpx)$. 

By flatness this is true for the morphism $\mathcal{X}\ra \cpx$, and blowing up $O(\Lambda)^-$ does not contribute new associated points, only the Cartier exceptional divisor $\pi^{-1}O(\Lambda)^-$. 

More precisely let $\mathcal{I}$ denote the ideal sheaf of $O(q)^- \subset \mathcal{X}$, and recall $\widehat{\mathcal{X}}$ is defined as $\textrm{Proj}\bigoplus_{d \geq 0}\mathcal{I}^d$. Any homogeneous zero divisor in the graded sheaf $\bigoplus_{d \geq 0}\mathcal{I}^d$ is already a zero divisor when regarded as an element of $\olo_\mathcal{X}$. On the other hand an associated point $\widehat{x} \in \widehat{\mathcal{X}}$ is by definition (following \cite{hart} III Corollary 9.6) a point for which every element of $\mathfrak{m}_{\widehat{x}}$ is a zero divisor. The natural map $\widehat{\mathcal{X}}\ra \mathcal{X}$ maps $\mathfrak{m}_{\widehat{x}}$ to its degree $0$ piece. Thus by the above remark $\widehat{x}$ necessarily maps to an associated point $x\in\mathcal{X}$. But $x$ maps to the generic point of $\textrm{Spec}(\cpx)$ by flatness, so the same is true for $\widehat{x}$.
\fine
\begin{remark} In both cases the assumption that $\alpha$ acts through $\Sl$ is not really restrictive. This can always be achieved by replacing $\mathcal{L}$ by some power and pulling back $\mathcal{X}$ by $z \mapsto z^d$ for some $d$. This gives a new test configuration for which $\alpha$ can be rescaled to act through $\Sl$ and for which the Futaki invariant is only multiplied by $d$, by Remark \ref{covering}.

This property of the Futaki invariant turns out to be important in our proof of Theorem \ref{main_teo}.
\end{remark}
\section{Proof of Theorem \ref{main_teo}}\label{proof} It will be enough to prove Proposition \ref{defo} and to apply the result of Arezzo and Pacard recalled as Theorem \ref{arezzo} below.

Thus let 
\[(\widehat{X}, L_{\gamma}) = (\Bl_q X, \pi^*L^{\gamma}\otimes\olo(-E)).\]
We need to show that $(\widehat{X}, L_{\gamma})$ is K-unstable for $\gamma \gg 0$. We will construct test configurations $(\mathcal{X}_{\gamma}, \mathcal{L}_{\gamma})$ for $(\widehat{X}, L_{\gamma})$ which have strictly negative Donaldson-Futaki invariant for $\gamma \gg 0$.\\

By assumption $(X, L)$ is properly semistable, so it admits a nontrivial test configuration $(\mathcal{X}, \mathcal{L})$ with $F(\mathcal{X}) = 0$.

Moreover we can assume that the induced $\cpx^*$-action on $H^0(\mathcal{X}_0, \mathcal{L}_0)^*$ is special linear. Indeed this can be achieved by taking some power $\mathcal{L}^r$ and a ramified cover $z \mapsto z^d$. The new Futaki invariant $F'$ still vanishes since $F' = d\cdot F = 0$.\\ 

We blow $\mathcal{X}$ up along the closure $O(q)^-$ of the orbit $O(q)$ of $q \in \mathcal{X}_1$ under the $\cpx^*$-action on $\mathcal{X}$, i.e. define
\begin{equation}
\mathcal{X}_{\gamma} = \widehat{\mathcal{X}} = \Bl_{O(q)^-}\mathcal{X}.
\end{equation}
Let $\olo_{\widehat{\mathcal{X}}}(1)$ denote the exceptional invertible sheaf on $\widehat{\mathcal{X}}$. We endow $\widehat{\mathcal{X}}$ with the polarisation 
\begin{equation}
\mathcal{L}_{\gamma} = \pi^*\mathcal{L}^{\gamma}\otimes\olo_{\widehat{\mathcal{X}}}(1). 
\end{equation}
Define the closed point $q_0 \in \mathcal{X}_0$ to be the specialisation
\[q_0 = \lim_{\lambda \ra 0}\lambda\cdot q.\]\\

Applying the blowup formula \ref{blowup} in this case gives
\begin{align*} 
F(\widehat{\mathcal{X}}_0, \pi^*\mathcal{L}^{\gamma}_0\otimes\olo_{\widehat{\mathcal{X}_0}}(1)) &= F(\mathcal{X}_0, \mathcal{L}_0) -\mathcal{CH}_{(\mathcal{X}_0, \mathcal{L}_0)}(q_0)\frac{\gamma^{1-n}}{2(n-2)!} + O(\gamma^{-n})\\
&= -\mathcal{CH}_{(\mathcal{X}_0, \mathcal{L}_0)}(q_0)\frac{\gamma^{1-n}}{2(n-2)!} + O(\gamma^{-n}). 
\end{align*}
In Proposition \ref{repulsive} below we will prove that for a very special $q \in \mathcal{X}_1 \cong X$, 
\[\mathcal{CH}_{(\mathcal{X}_0, \mathcal{L}_0)}(q_0) > 0.\]
This holds thanks to the assumption $F(\mathcal{X}) = 0$, or more generally $F(\mathcal{X}) \leq 0$. This is enough to settle Proposition \ref{defo}.\\

The final step for Theorem \ref{main_teo} is to show that the perturbation problem is unobstructed provided $\Aut(X, L)$ is discrete. This is precisely the content of a beautiful result of C. Arezzo and F. Pacard.
\begin{theorem}[Arezzo-Pacard \cite{arezzoI}]\label{arezzo} Let $(X, L)$ be a polarised manifold with a cscK metric in the class $c_1(L)$. Suppose $\Aut(X,L)$ is discrete and let $q \in X$ be any point. Then the blowup $\emph{Bl}_{q}X$ with exceptional divisor $E$ admits a cscK metric in the class $\gamma\pi^*c_1(L)-c_1(\olo(E))$ for $\gamma \gg 0$.
\end{theorem}
\begin{remark}
The Arezzo-Pacard theorem also holds in the K\"ahler case and, more importantly, even when $\mathfrak{aut}(X, L) \neq 0$, provided a suitable stability condition is satisfied. We refer to \cite{arezzoII}, \cite{io} for further discussion. 
\end{remark}
Thus the following Proposition will complete our proof(s). We believe it may also be of some independent interest.
\begin{proposition}\label{repulsive} Let $(\mathcal{X}, \mathcal{L})$ be a nonproduct test configuration for a polarised manifold $(X, L)$ with nonpositive Donaldson-Futaki invariant and suppose the induced $\cpx^*$-action on $H^0(\mathcal{X}_0, \mathcal{L}_0)^*$ is special linear. Then there exists $q \in \mathcal{X}_1\cong X$ such that $\mathcal{CH}_{(\mathcal{X}_0, \mathcal{L}_0)}(q_0) > 0$.
\end{proposition}
\dimo By the embedding Theorem \ref{embed} we reduce to the case of a nontrivial $\cpx^*$ acting on $\proj^N$ for some $N$, of the form $\textrm{Diag}(\lambda^{m_0}x_0, ... \lambda^{m_N}x_N)$, ordered by
\[m_0 \leq m_1 ...\,\leq m_N.\] 

Let $\{Z_i\}^k_{i = 1}$ be the distinct projective weight spaces, where $Z_{i}$ has weight $m_{i}$ (i.e. the induced action on $Z_i$ is trivial with weight $m_i$). Each $Z_i$ is a projective subspace of $\proj^N$, and the central fibre with its reduced induced structure $\mathcal{X}^{\textrm{red}}_0$ is a contained in $\textrm{Span}(Z_{i_1},...\,Z_{i_l})$ for some minimal flag $0 = i_1 < i_2 ... < i_l$.\\ 

\emph{The case $1 < l$.} In this case the induced action on closed points of $\mathcal{X}_0$ is nontrivial. Let $q \in \mathcal{X}_1$ be any point with 
\[\lim_{\lambda \ra 0} \lambda\cdot q = q_0 \in Z_{i_l}.\]
Such a point exists by minimality and because the specialisation of every point must lie in some $Z_j$. Since the action on $\mathcal{X}_0$ is induced from that on $\proj^N$, $q_0$ belongs to the totally repulsive fixed locus $R = \mathcal{X}_0 \cap Z_{i_l} \subset \mathcal{X}_0$. By this we mean that every closed point in $\mathcal{X}_0\setminus R$ specialises to a closed point in $\mathcal{X}_0\setminus R$. In particular the natural birational morphism $\mathcal{X}_0\dashrightarrow \textrm{Proj}(\bigoplus_{d}H^0(\mathcal{X}_0, \mathcal{L}^{\otimes d}_0)^{\cpx^*})$ blows up along $R$. So $q_0 \in R$ is an unstable point for the $\cpx^*$-action in the sense of geometric invariant theory. By the Hilbert-Mumford criterion the weight of the induced action on the line $\mathcal{L}_0|_{q_0}$ must be strictly positive. Since we are assuming that the induced action on $H^0(\mathcal{X}_0, \mathcal{L}_0)^*$ is special linear this weight coincides with the Chow weight, so $\mathcal{CH}_{(\mathcal{X}_0, \mathcal{L}_0)(q_0)} > 0$.\\ 
 
\emph{Degenerate case.} In the rest of the proof we will show that in the degenerate case $\mathcal{X}^{\textrm{red}}_0 \subset Z_0$ the Donaldson-Futaki invariant is strictly positive. Note that since by assumption the original $\cpx^*$-action on $\proj^N$ is nontrivial, $Z_0 \subset \proj^N$ is a proper projective subspace.\\

We digress for a moment to make the following observation: for any $\cpx^*$-action on $\proj^N$ with ordered weights $\{m_i\}$, and a smooth nondegenerate manifold $Y \subset \proj^N$, the map $\rho: Y \ni y \mapsto y_0 = \lim_{\lambda \ra 0} \lambda\cdot y$ is \emph{rational}, defined on the open dense set $\{y \in Y: \mu(y) = m_0\}$ of points with minimal Hilbert-Mumford weight. Indeed, in the above notation, generic points specialise to some point in the lowest fixed locus $Z_0$. In any case the map $\rho$ blows up exactly along loci where the Hilbert-Mumford weight jumps.\\

Going back to our discussion of the case $\mathcal{X}^{\textrm{red}}_0 \subset Z_0$, we see that this means precisely that all points of $\mathcal{X}_1$ have minimal Hilbert-Mumford weight $m_0$, so there is a well defined \emph{morphism}
\[\rho: \mathcal{X}_1 \ra Z_0.\] 
Moreover $\rho$ is a finite map: the pullback of $\mathcal{L}_0$ under $\rho$ is $L$ which is ample, therefore $\rho$ cannot contract a positive dimensional subscheme. If $\rho$ were an isomorphism on its image, it would fit in a $\cpx^*$-equivariant isomorphism $\mathcal{X} \cong X\times\cpx$. Therefore $\rho$ cannot be injective, either on closed points or tangent vectors. If, say, $\rho$ identifies distinct points $x_1, x_2$, this means that the $x_i$ specialise to the same $x$ under the $\cpx^*$-action; by flatness then the local ring $\olo_{\mathcal{X}_0,x }$ contains a nontrivial nilpotent pointing outwards of $Z_0$, i.e. the sheaf $\mathcal{I}_{\mathcal{X}_0 \cap Z_0}/\mathcal{I}_{\mathcal{X}_0}$ is nonzero. In other words $\mathcal{X}_0$ is not a closed subscheme of $Z_0$. The case when $\rho$ annihilates a tangent vector produces the same kind of nilpotent in the local ring of the limit, by specialisation.\\

To sum up, the central fibre $\mathcal{X}_0$ is nonreduced, containing nontrivial $Z_0$-orthogonal nilpotents. Equally important, the induced action on the closed subscheme $\mathcal{X}_0 \cap Z_0 \subset \mathcal{X}_0$ is trivial. The proof will be completed by a weight computation.\\ 

\emph{Donaldson-Futaki invariant.} Suppose $Z_0 \subset \proj^N$ has projective coordinates $[x_1 : ... : x_r]$, i.e. it is cut out by $\{x_{r+1} = ... = x_N = 0\}$. We change the linearisation by changing the representation of the $\cpx^*$-action, to make it of the form
\begin{equation}\label{linearisation}
[x_0: ... x_r : x_{r+1}: ...: x_N] \mapsto [x_0: ... x_r : \lambda^{m_{r+1}-m_0}x_{r+1}: ...: \lambda^{m_N - m_0}x_N],
\end{equation}
and recall $m_{r+i} > m_0$ for all $i > 0$. It is possible that the induced action on $H^0(\mathcal{X}_0, \mathcal{L}_0)^*$ will not be special linear anymore, however this does not affect the Donaldson-Futaki invariant.

Note that for all large $k$,  
\begin{equation}\label{sections}
H^0(\proj^N, \olo(k)) \ra H^0(\mathcal{X}_0, \mathcal{L}^k_0)\ra H^1(\mathcal{I}_{\mathcal{X}_0}(k)) = 0.
\end{equation}
By \ref{sections}, our geometric description of $\mathcal{X}_0$ and the choice of linearisation \ref{linearisation} we see that any section $\xi \in H^0(\mathcal{X}_0, \mathcal{L}^k_0)$ has nonnegative weight under the induced $\cpx^*$-action. The section $\xi$ can only have strictly positive weight if it is of the form $x_{r+i}\cdot f$ for some $i > 0$. Moreover we know there exists an integer $a > 0$ such that $x^a_{r + i}|_{\mathcal{X}_0} = 0$ for all $i > 0$. Let $w(k)$ denote the total weight of the action on $H^0(\mathcal{X}_0, \mathcal{L}^k_0)$, i.e. the induced weight on the line $\Lambda^{P(k)}H^0(\mathcal{X}_0, \mathcal{L}^k_0)$, where $P(k) = h^0(\mathcal{X}_0, \mathcal{L}^k_0)$ is the Hilbert polynomial. Our discussion implies the upper bound
\begin{equation}
w(k)\leq C (P(k-1) + ...+ P(k-a))
\end{equation}
for some $C > 0$, independent of $k$. In particular, 
\begin{equation}\label{lo_bound}
w(k) = O(k^n).
\end{equation}
On the other hand we can look at just one section $x_{r + i}$, $i > 0$ with $x_{r + i}|_{\mathcal{X}_0} \neq 0$. This gives a lower bound
\begin{equation}
w(k) \geq C\cdot P(k-1)
\end{equation}
for some $C > 0$, independent of $k$. So we see that
\begin{equation}
\frac{w(k)}{k P(k)} \geq \frac{C'}{k}.
\end{equation}
holds for $k \gg 0$ and some $C' > 0$ independent of $k$. Together with 
\begin{equation}
\frac{w(k)}{k P(k)} = O(k^{-1})
\end{equation}
which follows from \ref{lo_bound} this implies 
\begin{equation}
\frac{w(k)}{k P(k)} = \frac{C''}{k} + O(k^{-2})
\end{equation}
for some $C'' > 0$ independent of $k$.

By definition of Donaldson-Futaki invariant, this immediately implies
\[F(\mathcal{X}) \geq C'' > 0, \]
a contradiction.
\fine
\begin{remark}\label{flags} One can characterise the degenerate case in the above proof  more precisely. 

As observed by Ross-Thomas \cite{ross_thomas} Section 3 a result of Mumford implies that any test configuration $(\mathcal{X}, \mathcal{L})$ for $(X, L)$ is a contraction of some blowup of $X \times \cpx$ in a flag of $\cpx^*$-invariant closed subschemes supported in some thickening of $X \times \{0\}$.

The existence of the map $\rho: \mathcal{X}_1 \ra Z_0$ means precisely that in this Mumford representation of $\mathcal{X}$ no blowup occurs, i.e. $\mathcal{X}$ is a contraction of the product $X \times \cpx$.

Define a map $\nu: X \times \cpx \ra \mathcal{X}$ by $\nu(x, \lambda) = \lambda\cdot x$ away from $X \times \{0\}$, $\nu = \rho$ on $X \times \{0\}$. This is a well defined morphism, and since $\rho$ is finite, $\nu$ is precisely the \emph{normalisation} of $\mathcal{X}$.

So in the degenerate case $\mathcal{X}^{\textrm{red}}_0 \subset Z_0$ the normalisation of $\mathcal{X}$ is $X \times \cpx$.

Ross-Thomas \cite{ross_thomas} Proposition 5.1 proved the general result that normalising a test configuration reduces the Donaldson-Futaki invariant. This already implies $F \geq 0$ in the degenerate case, since the induced action on $X \times \cpx$ must have vanishing Futaki invariant. In our special case our direct proof yields the strict inequality we need.
\end{remark}
\begin{remark} The result of Mumford mentioned above states more precisely that any test configuration $(\mathcal{X}, \mathcal{L})$ for $(X, L)$ is a contraction of the blowup of $X \times \cpx$ in an ideal sheaf
\[\mathbf{I}_r = \mathcal{I}_0 + t\mathcal{I}_1 + ... + t^{r-1}\mathcal{I}_{r-1} + (t^r)\]
where $\mathcal{I}_0 \subseteq ... \subseteq \mathcal{I}_{r-1} \subset \olo_X$ correspond to a descending flag of closed subschemes $Z_0 \supseteq ... \supseteq Z_{r-1}$. The action on $(\mathcal{X}, \mathcal{L})$ is the natural one induced from the trivial action on $X \times \cpx$. 

Suppose now that $F(\mathcal{X}) = 0$ and that \emph{no contraction} occurs in Mumford's representation. 

Then in Proposition \ref{repulsive} we can simply choose any closed point $q \in Z_{r-1}$. This is because the proper transform of $Z_{r-1}\times\cpx$ cuts $\mathcal{X}_0$ in the totally repulsive locus for the induced action, i.e. the action flows every closed point in $\mathcal{X}_0$ outside this locus to the proper transform of $X \times \{0\}$.

Conversely blowing up $q \in X \setminus Z_0$ only increases the Donaldson-Futaki invariant (at least asymptotically).

For example K-stability with respect to test configurations with $r = 1$ and no contraction is known as Ross-Thomas \emph{slope stability} \cite{ross_thomas} and has found interesting applications to cscK metrics. In particular this discussion gives a simpler proof that a cscK polarised manifold with discrete automorphisms is slope stable.
\end{remark}
\begin{remark}\label{uniform}
A refinement of Conjecture \ref{don_conj} was proposed by G. Sz\'ekelyhidi. If $\om \in c_1(L)$ is cscK there should be a \emph{strictly positive lower bound} for a suitable normalisation of $F$ over all nonproduct test configurations. This condition is called \emph{uniform K-polystability}. In \cite{gabor} Section 3.1.1 it is shown that the correct normalisation in the case of algebraic surfaces coincides with that of Theorem \ref{calabi}, namely $\frac{F(\mathcal{X})}{\|\mathcal{X}\|}$. For toric surfaces K-polystability implies uniform K-polystability with respect to toric test configurations; this is shown in \cite{gabor} Section 4.2. It seems clear however that the proof presented here cannot be refined to yield uniform K-stability for surfaces.
\end{remark}
\begin{example}[Del Pezzo surfaces] Del Pezzo surfaces played an important role in the development of the subject. By the work of Tian and others all smooth Del Pezzo surfaces $V_d$ of degree $d \leq 6$ admit a K\"ahler-Einstein metric. For $d \leq 5$, $V_d$ has discrete automorphism group. K-stability in the sense of Tian for $V_d$, $d \leq 5$ follows from \cite{tian} Theorem 1.2. K-stability with respect to ``good" test configurations follows from \cite{paul} Theorem 2. 

Our Theorem \ref{main_teo} refines this to K-stability in the sense of Donaldson. 

Moreover Theorem \ref{main_teo} also applies to polarisations on $V_d$, $d \leq 5$ for which the exceptional divisors have sufficiently small volume, thanks to the results of Arezzo and Pacard \cite{arezzoII}. 
\end{example}

\vskip.3cm

\noindent Universit\`a di Pavia, Via Ferrata 1 27100 Pavia, ITALY\\
and\\
Imperial College, London SW7 2AZ, UK.\\
\textit{E-mail:} jacopo.stoppa\texttt{@unipv.it}
\end{document}